\documentclass[11pt]{amsart}
\usepackage{amssymb}
\usepackage{amsmath}
\usepackage[active]{srcltx}
\usepackage{t1enc}
\usepackage[latin2]{inputenc}
\usepackage{verbatim}
\usepackage{amsmath,amsfonts,amssymb,amsthm}
\usepackage[mathcal]{eucal}
\usepackage{enumerate}
\usepackage[centertags]{amsmath}
\usepackage{graphics}

\newtheorem{theorem}{Theorem}

\newtheorem{lemma}{Lemma}

\begin{document}
\author{Istv\'an Blahota, Giorgi Tephnadze and Rodolfo Toledo}
\title[Ces\`aro means]{Strong convergence theorem of Ces\`aro means with
respect to the Walsh system}
\address{I. Blahota, Institute of Mathematics and Computer Sciences, College
of Ny\'\i regyh\'aza, P.O. Box 166, Ny\'\i regyh\'aza, H-4400, Hungary}
\email{blahota@nyf.hu}
\address{G. Tephnadze, Department of Mathematics, Faculty of Exact and
Natural Sciences, Tbilisi State University, Chavchavadze str. 1, Tbilisi
0128, Georgia and Department of Engineering Sciences and Mathematics, Lule%
\aa\ University of Technology, SE-971 87 Lule\aa , Sweden}
\email{giorgitephnadze@gmail.com}
\address{R. Toledo, Institute of Mathematics and Computer Sciences, College
of Ny\'\i regyh\'aza, P.O. Box 166, Ny\'\i regyh\'aza, H-4400, Hungary}
\email{toledo@nyf.hu}
\thanks{The research was supported by project
T\'AMOP-4.2.2.A-11/1/KONV-2012-0051 and by Shota Rustaveli National Science
Foundation grant no.52/54 (Bounded operators on the martingale Hardy
spaces). }
\date{}

\begin{abstract}
We prove that Ces\`{a}ro means of one-dimensional Walsh-Fourier series are
uniformly bounded operators in the martingale Hardy space $H_{p}$ for $%
0<p<1/\left( 1+\alpha \right).$
\end{abstract}

\maketitle
\date{}

\section{Introduction}

\ The definitions and notations used in this introduction can be found in
the next Section. \ It is well-known (see, e.g., \cite[p.125]{G-E-S}) that
Walsh-Paley system is not a Schauder basis in the space $L_{1}\left(
G\right) $. Moreover, there is a function $F$ in the dyadic Hardy space $%
H_{1}\left( G\right) $, such that the partial sums of the Walsh-Fourier
series of $F$ are not bounded in the $L_{1}$-norm. However, in Simon \cite%
{Si3} the following estimation was obtained: for all $F\in H_{1}\left(
G\right)$ 
\begin{equation*}
\frac{1}{\log n}\overset{n}{\underset{k=1}{\sum }}\frac{\left\Vert
S_{k}F\right\Vert _{1}}{k}\leq c\left\Vert F\right\Vert _{H_{1}},\text{ \ \
\ }\left( n=2,3,\ldots \right),
\end{equation*}%
where $S_{k}F$ denotes the $k$-th partial sum of the Walsh-Fourier series of 
$F$ (For the trigonometric analogue see in Smith \cite{sm}, for the Vilenkin
system in G\'at \cite{gat1}, for a more general, so-called Vilenkin-like
system in Blahota \cite{b}). Simon \cite{si1} (see also \cite{tep6} and \cite%
{We}) proved that there exists an absolute constant $c_{p},$ depending only
on $p,$ such that%
\begin{equation}
\frac{1}{\log ^{\left[ p\right] }n}\overset{n}{\underset{k=1}{\sum }}\frac{%
\left\Vert S_{k}F\right\Vert _{p}^{p}}{k^{2-p}}\leq c_{p}\left\Vert
F\right\Vert _{H_{p}}^{p},\ \ \left( 0<p\leq 1,\text{ }n=2,3,\ldots \right) ,
\label{1ccc}
\end{equation}%
for all $F\in H_{p}$, where $\left[ p\right] $ denotes integer part of $p.$

In \cite{tep4} it was proven that sequence $\left\{ 1/k^{2-p}\right\}
_{k=1}^{\infty }$ $\left( 0<p<1\right) $ \ in (\ref{1ccc}) is given exactly.

Weisz \cite{We3} considered the norm convergence of Fej\'er means of
Walsh-Fourier series and proved that 
\begin{equation}
\left\Vert \sigma _{n}F\right\Vert _{H_{p}}\leq c_{p}\left\Vert F\right\Vert
_{H_{p}},\text{ \ \ }F\in H_{p},\text{\ \ \ \ }\left( 1/2<p<\infty ,\text{ \ 
}n=1,2,3,\ldots \right) ,  \label{1ac}
\end{equation}%
where the constant $c_{p}>0$ depends only on $p$.

Inequality (\ref{1ac}) immediately implies that 
\begin{equation*}
\frac{1}{n^{2p-1}}\overset{n}{\underset{k=1}{\sum }}\frac{\left\Vert \sigma
_{k}F\right\Vert _{H_{p}}^{p}}{k^{2-2p}}\leq c_{p}\left\Vert F\right\Vert
_{H_{p}}^{p},\text{ \ \ \ }\left( 1/2<p<\infty \right) .
\end{equation*}

If (\ref{1ac}) also hold, for $0<p\leq 1/2,$ then we would have 
\begin{equation}
\frac{1}{\log ^{\left[ 1/2+p\right] }n}\overset{n}{\underset{k=1}{\sum }}%
\frac{\left\Vert \sigma _{k}F\right\Vert _{H_{p}}^{p}}{k^{2-2p}}\leq
c_{p}\left\Vert F\right\Vert _{H_{p}}^{p},\text{ \ \ \ }\left( 0<p\leq 1/2,%
\text{ }n=2,3,\ldots \right) ,  \label{2cc}
\end{equation}%
but in \cite{tep1} it was proven that the assumption $p>1/2$ is essential.
In particular, there was proven that there exists a martingale $F\in H_{p}$ $%
\left( 0<p\leq 1/2\right) ,$ such that $\sup_{n}\left\Vert \sigma
_{n}F\right\Vert _{p}=+\infty .$

However, in \cite{tep5} (see also \cite{bt2}) it was proven that (\ref{2cc})
holds, though (\ref{1ac}) is not true for $0<p\leq 1/2.$

The weak-type (1,1) inequality for the maximal operator of Fej\'er means $%
\sigma ^{\ast }$ can be found in Schipp \cite{Sc} (see also \cite{PS}).
Fujji \cite{Fu} and Simon \cite{Si2} verified that $\sigma ^{\ast }$ is
bounded from $H_{1}$ to $L_{1}$. Weisz \cite{We2} generalized this result
and proved the boundedness of $\sigma ^{\ast }$ from the space $H_{p}$ to
the space $L_{p}$ for $p>1/2$. Simon \cite{Si1} gave a counterexample, which
shows that boundedness does not hold for $0<p<1/2.$ The counterexample for $%
p=1/2$ is due to Goginava \cite{GoAMH} (see also \cite{BGG2}). Weisz \cite%
{we4} proved that $\sigma ^{\ast }$ is bounded from the Hardy space $H_{1/2}$
to the space $L_{1/2,\infty }$. In \cite{tep2, tep3} it was proven that the
maximal operators $\widetilde{\sigma }_{p}^{\ast }$ defined by 
\begin{equation}
\widetilde{\sigma }_{p}^{\ast }F:=\sup_{n\in \mathbb{N}}\frac{\left\vert
\sigma _{n}F\right\vert }{n^{1/p-2}\log ^{2\left[ 1/2+p\right] }n},\text{ \ }%
\left( 0<p\leq 1/2,\text{ }n=2,3,\ldots \right)  \label{1}
\end{equation}%
is bounded from the Hardy space $H_{p}$ to the space $L_{p},$ where $F\in
H_{p}$ and $\left[ 1/2+p\right] $ denotes integer part of $1/2+p.$ Moreover,
there was also shown that sequence $\left\{ n^{1/p-2}\log ^{2\left[ 1/2+p%
\right] }n\text{ \ }:\text{ \ }n=2,3,\ldots \right\} $ in (\ref{1}) can not
be improved.

Weisz \cite{we6} proved that the maximal operator $\sigma ^{\alpha ,\ast }$ $%
\left( 0<\alpha <1\right) $ of the Ces\`{a}ro means of Walsh system is
bounded from the martingale space $H_{p}$ to the space $L_{p}$ for $%
p>1/\left( 1+\alpha \right) .$ Goginava \cite{gog4} gave a counterexample,
which shows that the boundedness does not hold for $0<p\leq 1/\left(
1+\alpha \right) .$ Recently, Weisz and Simon \cite{sw} show that the
maximal operator $\sigma ^{\alpha ,\ast }$ is bounded from the Hardy space $%
H_{1/\left( 1+\alpha \right) }$ to the space $L_{1/\left( 1+\alpha \right)
,\infty }$. An analogical result for Walsh-Kaczmarz system was proven in 
\cite{gg}.

In \cite{gog8} Goginava investigated the behaviour of Ces\`aro means of
Walsh-Fourier series in detail. For some approximation properties of the two
dimensional case see paper of Nagy \cite{na}.

The main aim of this paper is to generalize estimate (\ref{2cc}) for Ces\`{a}%
ro means, when $0<p<1/\left( 1+\alpha \right) .$ We also consider the
weighted maximal operator of $\left( C,\alpha \right) $ means and proved
some new $\left( H_{p},L_{p}\right) $-type inequalities for it.

We note that the case $p=1/\left( 1+\alpha \right) $ was considered in \cite%
{bt1}.

\section{Definitions and Notations}

Let $\mathbb{N}_{+}$ denote the set of the positive integers, $\mathbb{N}:=%
\mathbb{N}_{+}\cup \{0\}.$ Denote by $\mathbb{Z}_{2}$ the discrete cyclic
group of order 2, that is $\mathbb{Z}_{2}:=\{0,1\},$ where the group
operation is the modulo 2 addition and every subset is open. The Haar
measure on $\mathbb{Z}_{2}$ is given so that the measure of a singleton is $%
1/2$.

Define the group $G$ as the complete direct product of the group $\mathbb{Z}%
_{2}$ with the product of the discrete topologies of $\mathbb{Z}_{2}$'s. The
elements of $G$ are represented by sequences 
\begin{equation*}
x:=(x_{0},x_{1},\ldots ,x_{k},\ldots )\qquad \left( x_{k}=0,1\right) .
\end{equation*}
It is easy to give a base for the neighborhood of $G$ 
\begin{eqnarray*}
I_{0}\left( x\right) &:&=G, \\
I_{n}(x) &:&=\{y\in G\mid y_{0}=x_{0},\ldots ,y_{n-1}=x_{n-1}\}\ (x\in G,\
n\in \mathbb{N}).
\end{eqnarray*}%
Denote $I_{n}:=I_{n}\left( 0\right) \ $and $\overline{I_{n}}:=G$ $\backslash 
$ $I_{n}$. Let 
\begin{equation*}
e_{n}:=\left( 0,\dots,0,x_{n}=1,0,\dots\right) \in G\qquad \left( n\in 
\mathbb{N}\right)
\end{equation*}

Denote 
\begin{equation*}
I_{M}^{k,l}:=\left\{ 
\begin{array}{ll}
I_{M}(0,\ldots ,0,x_{k}=1,0,\dots ,0,x_{l}=1,x_{l+1},\ldots ,x_{M-1}), & 
k<l<M, \\ 
I_{M}(0,\ldots ,0,x_{k}=1,0,\ldots ,0), & l=M.%
\end{array}%
\right.
\end{equation*}
It is evident 
\begin{equation}
\overline{I_{M}}=\left( \overset{M-2}{\underset{k=0}{\bigcup }}\overset{M-1}{%
\underset{l=k+1}{\bigcup }}I_{M}^{k,l}\right) \bigcup \left( \underset{k=0}{%
\bigcup\limits^{M-1}}I_{M}^{k,M}\right) .  \label{2}
\end{equation}

If $n\in \mathbb{N},$ then every $n$ can be uniquely expressed as $%
n=\sum_{j=0}^{\infty }n_{j}2^{j},$ where $n_{j}\in Z_{2}$ $(j\in \mathbb{N})$
and only finite number of $n_{j}$'s differ from zero, that is, $n$ is
expressed in the number system of base 2. Let $\left\vert n\right\vert
:=\max $ $\{j\in \mathbb{N},$ $n_{j}\neq 0\},$ that is $2^{\left\vert
n\right\vert }\leq n\leq 2^{\left\vert n\right\vert +1}.$

The norm (or quasi-norm) of the space $L_{p}(G)$ is defined by 
\begin{equation*}
\left\Vert f\right\Vert _{p}:=\left( \int_{G}\left\vert f\right\vert
^{p}d\mu \right) ^{1/p},\ \left( 0<p<\infty \right) .
\end{equation*}
The space $L_{p,\infty }\left( G\right) $ consists of all measurable
functions $f$, for which

\begin{equation*}
\left\Vert f\right\Vert _{L_{p,\infty }(G)}:=\underset{\lambda >0}{\sup }%
\lambda \mu \left( f>\lambda \right) ^{1/p}<\infty .
\end{equation*}

Next, we introduce on $G$ an orthonormal system which is called the \textit{%
Walsh system}. At first, define the functions $r_{k}\left( x\right)
:G\rightarrow \mathbb{C} $, the so-called Rademacher functions as 
\begin{equation*}
r_{k}\left( x\right) :=\left( -1\right) ^{x_{k}}\qquad \left( x\in G,\ k\in 
\mathbb{N}\right) .
\end{equation*}

Now, define the Walsh system $w:=(w_{n}:n\in \mathbb{N})$ on $G$ as: 
\begin{equation*}
w_{n}(x):=\prod_{k=0}^{\infty}r_{k}^{n_{k}}\left( x\right) =r_{\left\vert
n\right\vert }\left( x\right) \left( -1\right) ^{\underset{k=0}{\overset{%
\left\vert n\right\vert -1}{\sum }}n_{k}x_{k}}\qquad \left( n\in \mathbb{N}%
\right) .
\end{equation*}
The Walsh system is orthonormal and complete in $L_{2}\left( G\right) $
(see, e.g., \cite{Vi}).

If $\ f\in L_{1}\left( G\right) ,$ then we can establish Fourier
coefficients, partial sums of the Fourier series, Fej\'er means, Dirichlet
and Fej\'er kernels in the usual manner: 
\begin{eqnarray*}
\widehat{f}\left( n\right) &:&=\int_{G}fw_{n}d\mu \,,\,\,\ \ \ \ \ \ \ \
\left( n\in \mathbb{N}\right) , \\
S_{n}f &:&=\sum_{k=0}^{n-1}\widehat{f}\left( k\right) w_{k},\text{ \ \ \ \ }%
\left( n\in \mathbb{N}_{+}\right) , \\
\sigma _{n}f &:&=\frac{1}{n}\sum_{k=1}^{n}S_{k}f,\ \ \ \ \ \ \ \ \left( n\in 
\mathbb{N}_{+}\right) , \\
D_{n} &:&=\sum_{k=0}^{n-1}w_{k},\,\ \ \ \ \ \ \ \ \ \ \text{\ \ }\left( n\in 
\mathbb{N}_{+}\right) , \\
K_{n} &:&=\frac{1}{n}\overset{n}{\underset{k=1}{\sum }}D_{k},\ \ \ \ \ \
\,\left( n\in \mathbb{N}_{+}\right) ,
\end{eqnarray*}%
respectively. Recall that (see e.g., \cite{sws}) 
\begin{equation}
D_{2^{n}}\left( x\right) =\left\{ 
\begin{array}{ll}
2^{n}, & \text{ if }x\in I_{n}, \\ 
0, & \text{ if }x\notin I_{n}.%
\end{array}%
\right.  \label{1dn}
\end{equation}

The Ces\`{a}ro means ($\left( C,\alpha \right) $-means) are defined as 
\begin{equation*}
\sigma _{n}^{\alpha }f:=\frac{1}{A_{n}^{\alpha }}\overset{n}{\underset{k=1}{%
\sum }}A_{n-k}^{\alpha -1}S_{k}f,
\end{equation*}%
where 
\begin{equation}
A_{0}^{\alpha }:=1,\qquad A_{n}^{\alpha }:=\frac{\left( \alpha +1\right)
\dots \left( \alpha +n\right) }{n!}\ \ \alpha \neq -1,-2,\dots  \label{1d}
\end{equation}
It is well known that 
\begin{equation*}
A_{n}^{\alpha }=\sum_{k=0}^{n}A_{n-k}^{\alpha -1},\ A_{n}^{\alpha
}-A_{n-1}^{\alpha }=A_{n}^{\alpha -1},\ A_{n}^{\alpha }\backsim n^{\alpha }
\end{equation*}%
and 
\begin{equation}
\sup_{n}\int_{G}\left\vert K_{n}^{\alpha }\right\vert d\mu \leq c<\infty ,
\label{4}
\end{equation}%
where $K_{n}^{\alpha }$ is $n$ th Ces\`{a}ro kernel.

The $\sigma $-algebra generated by the intervals $\{I_{n}(x):x\in G\}$ will
be denoted by $F_{n}\ (n\in \mathbb{N})$. Denote by $F=(F_{n},n\in\mathbb{N})
$ the martingale with respect to $F_{n}\ (n\in\mathbb{N}) $ (for details
see, e.g., \cite{We1}).

The maximal function of a martingale $F$ is defined by

\begin{equation*}
F^{\ast }:=\sup_{n\in \mathbb{N}}\left\vert F_{n}\right\vert .
\end{equation*}
In the case $f\in L_{1}\left( G\right) ,$ the maximal functions are also be
given by

\begin{equation*}
f^{\ast }\left( x\right) =\sup\limits_{n\in \mathbb{N}}\frac{1}{\mu \left(
I_{n}\left( x\right) \right) }\left\vert \int\limits_{I_{n}\left( x\right)
}f\left( u\right) d\mu \left( u\right) \right\vert .
\end{equation*}
For $0<p<\infty ,$ the Hardy martingale spaces $H_{p}\left( G\right) $
consist of all martingales such that

\begin{equation*}
\left\Vert F\right\Vert _{H_{p}}:=\left\Vert F^{\ast }\right\Vert
_{p}<\infty .
\end{equation*}

A bounded measurable function $a$ is a $p$-atom, if there exist a dyadic
interval $I$ such that \qquad 
\begin{equation*}
\int_{I}ad\mu =0,\text{ \ \ }\left\Vert a\right\Vert _{\infty }\leq \mu
\left( I\right) ^{-1/p},\text{ \ \ \ supp}\left( a\right) \subset I.
\end{equation*}

It is easy to check that for every martingale $F=\left( F_{n},n\in \mathbb{N}%
\right) $ and for every $k\in \mathbb{N}$ the limit 
\begin{equation}
\widehat{F}\left( k\right) :=\lim_{n\rightarrow \infty
}\int_{G}F_{n}w_{k}d\mu  \label{3a}
\end{equation}%
exists and it is called the $k$-th Walsh-Fourier coefficients of $F.$

Denote by ${\mathcal{A}}_{n}$ the $\sigma$-algebra generated by the sets $%
I_{n}(x)\,(x\in G,\ n\in {\mathbb{N}})$. If $F:=$ $\left( S_{2^{n}}f:n\in 
\mathbb{N}\right) $ is the regular martingale generated by $f\in L_{1}\left(
G\right)$, then 
\begin{equation*}
\widehat{F}\left( k\right) =\int_{G}fw_{k}d\mu =:\widehat{f}\left( k\right)
,\qquad k\in \mathbb{N}.
\end{equation*}

For $0<\alpha \leq 1,$ let consider maximal operators 
\begin{equation*}
\sigma ^{\alpha ,\ast }F:=\sup_{n\in \mathbb{N}}\left\vert \sigma
_{n}^{\alpha }F\right\vert ,\text{ \ }\overset{\sim }{\sigma }_{p}^{\alpha
,\ast }F:=\sup_{n\in \mathbb{N}}\frac{\left\vert \sigma _{n}^{\alpha
}F\right\vert }{\left( n+1\right) ^{1/p-1-\alpha }},\ 0<p<1/\left( 1+\alpha
\right) .
\end{equation*}

For the martingale 
\begin{equation*}
F=\sum_{n=0}^{\infty }\left( F_{n}-F_{n-1}\right)
\end{equation*}%
the conjugate transforms are defined as 
\begin{equation*}
\widetilde{F^{\left( t\right) }}:=\sum_{n=0}^{\infty }r_{n}\left( t\right)
\left( F_{n}-F_{n-1}\right) ,
\end{equation*}%
where $t\in G$ is fixed. Note that $\widetilde{F^{\left( 0\right) }}=F.$

As it is well-known (see, e.g., \cite{We1}) 
\begin{equation}
\left\Vert \widetilde{F^{\left( t\right) }}\right\Vert _{H_{p}}=\left\Vert
F\right\Vert _{H_{p}},\text{ \ \ }\left\Vert F\right\Vert _{H_{p}}^{p}\sim
\int_{G}\left\Vert \widetilde{F^{\left( t\right) }}\right\Vert _{p}^{p}dt,%
\text{ \ \ \ }\widetilde{\sigma _{m}^{\alpha }F^{\left( t\right) }}=\sigma
_{m}^{\alpha }\widetilde{F^{\left( t\right) }}.  \label{5.1}
\end{equation}

\section{Formulation of main results}

\begin{theorem}
\textbf{\label{Theorem1}}a) Let $0<\alpha <1$ and $0<p<1/(1+\alpha )$. Then
there exists absolute constant $c_{\alpha ,p}$, depending on $\alpha $ and $%
p $, such that for all $F\in H_{p}(G)$ 
\begin{equation*}
\left\Vert \overset{\sim }{\sigma }_{p}^{\alpha ,\ast }F\right\Vert _{p}\leq
c_{\alpha ,p}\left\Vert F\right\Vert _{H_{p}}.
\end{equation*}

b) Let $0<\alpha <1$, $0$\thinspace $<p\,<1/\left( 1+\alpha \right) $ and $%
\varphi :\mathbb{N}_{+}\rightarrow \lbrack 1,$ $\infty )$ be a nondecreasing
function satisfying the condition

\begin{equation}
\overline{\lim_{n\rightarrow \infty }}\frac{n^{1/p-1-\alpha }}{\varphi
\left( n\right) }=\infty .  \label{6}
\end{equation}
Then the maximal operator 
\begin{equation*}
\sup_{n\in \mathbb{N}}\frac{|\sigma_{n}^{\alpha}f|}{\varphi(n)}
\end{equation*}
is not bounded from the Hardy space $H_{p}(G)$ to the space $L_{p}(G)$.
\end{theorem}

\begin{theorem}
\textbf{\label{Theorem2}}Let $0<\alpha <1$ and $0<p<1/(1+\alpha )$. Then
there exists absolute constant $c_{\alpha ,p}$, depending on $\alpha $ and $%
p $, such that for all $F\in H_{p}$ 
\begin{equation*}
\overset{\infty }{\underset{m=1}{\sum }}\frac{\left\Vert \sigma _{m}^{\alpha
}F\right\Vert _{H_{p}}^{p}}{m^{2-\left( 1+\alpha \right) p}}\leq c_{\alpha
,p}\left\Vert F\right\Vert _{H_{p}}^{p}.
\end{equation*}
\end{theorem}

\section{Auxiliary Propositions}

The dyadic Hardy martingale spaces $H_{p}\left( G\right) $ have an atomic
characterization, when $0<p\leq 1$:

\begin{lemma}
(\textrm{Weisz\label{W} \cite{We5}}) A martingale $F=\left( F_{n},\ n\in 
\mathbb{N}\right) $ is in $H_{p}$ $\left( 0<p\leq 1\right) $ if and only if
there exists a sequence $\left( a_{k},\ k\in \mathbb{N}\right) $ of $p$%
-atoms and a sequence $\left( \mu _{k},\ k\in \mathbb{N}\right) $ of a real
numbers, such that for every $n\in \mathbb{N}$

\begin{equation}
\qquad \sum_{k=0}^{\infty }\mu _{k}S_{2^{n}}a_{k}=F_{n},  \label{2A}
\end{equation}%
\begin{equation*}
\text{\ }\sum_{k=0}^{\infty }\left\vert \mu _{k}\right\vert ^{p}<\infty .
\end{equation*}
Moreover, 
\begin{equation*}
\left\Vert F\right\Vert _{H_{p}}\backsim \inf \left( \sum_{k=0}^{\infty
}\left\vert \mu _{k}\right\vert ^{p}\right) ^{1/p},
\end{equation*}%
where the infimum is taken over all decompositions of $F$ of the form (\ref%
{2A}).
\end{lemma}

By using Lemma \ref{W} we can easily proved the following:

\begin{lemma}
(\textrm{Weisz\label{W1}} \cite{We1}) Suppose that an operator $T$ is $%
\sigma $-linear and for some $0<p\leq 1$ 
\begin{equation*}
\int\limits_{\overset{-}{I}}\left\vert Ta\right\vert ^{p}d\mu \leq
c_{p}<\infty ,
\end{equation*}%
for every $p$-atom $a$, where $I$ denote the support of the atom. If $T$ is
bounded from $L_{\infty }$ to $L_{\infty },$ then 
\begin{equation*}
\left\Vert Tf\right\Vert _{p}\leq c_{p}\left\Vert f\right\Vert _{H_{p}}.
\end{equation*}
\end{lemma}

To prove our main results we also need the following estimations:

\begin{lemma}
\textbf{\label{bt1}}\cite{bt1} Let $0<\alpha <1$ and $n>2^{M}.$ Then 
\begin{equation*}
\int_{I_{M}}\left\vert K_{n}^{\alpha }\left( x+t\right) \right\vert d\mu
\left( t\right) \leq \frac{c_{\alpha }2^{\alpha l+k}}{n^{\alpha }2^{M}},
\end{equation*}%
$\text{ }$for $x\in I_{l+1}\left( e_{k}+e_{l}\right) ,$ $(k=0,\ldots
,M-2,l=k+1,\ldots ,M-1)$ and 
\begin{equation*}
\int_{I_{M}}\left\vert K_{n}^{\alpha }\left( x+t\right) \right\vert d\mu
\left( t\right) \leq \frac{c_{\alpha }2^{k}}{2^{M}},
\end{equation*}%
\textit{for} \ $x\in I_{M}\left( e_{k}\right) ,$ $(k=0,\ldots ,M-1).$
\end{lemma}

\section{Proof of Theorems}

\textbf{Proof of Theorem 1.} Since $\sigma _{n}$ is bounded from $L_{\infty
} $ to $L_{\infty }$ (the boundedness follows from (\ref{4})) according to
Lemma \ref{W1} the proof of Theorem \ref{Theorem1} will be complete if we
show 
\begin{equation*}
\sup \int_{\overline{I_{M}}}\left\vert \overset{\sim }{\sigma }_{p}^{\alpha
,\ast }a\right\vert ^{p}d\mu <\infty ,
\end{equation*}%
where the supremum is taken over all $p$-atoms $a$. We may assume that $a$
is an arbitrary $p$-atom, with support$\ I,$ $\mu \left( I\right) =2^{-M}$
and $I=I_{M}.$ It is easy to see that $\sigma _{n}^{\alpha }\left( a\right)
=0,$ when $n\leq 2^{M}.$ Therefore, we can suppose that $n>2^{M}.$

Let $x\in I_{M}.$ Since $\left\Vert a\right\Vert _{\infty }\leq c2^{M/p}$ we
obtain 
\begin{equation*}
\left\vert \sigma _{n}^{\alpha }a\left( x\right) \right\vert \leq
\int_{I_{M}}\left\vert a\left( t\right) \right\vert \left\vert K_{n}^{\alpha
}\left( x+t\right) \right\vert d\mu \left( t\right)
\end{equation*}%
\begin{equation*}
\leq \left\Vert a\left( x\right) \right\Vert _{\infty
}\int_{I_{M}}\left\vert K_{n}^{\alpha }\left( x+t\right) \right\vert d\mu
\left( t\right)
\end{equation*}%
\begin{equation*}
\leq c_{\alpha }2^{M/p}\int_{I_{M}}\left\vert K_{n}^{\alpha }\left(
x+t\right) \right\vert d\mu \left( t\right) .
\end{equation*}

Let $x\in I_{M}^{k,l},\,0\leq k<l<M.$ Then from Lemma \ref{bt1} we get 
\begin{equation}
\left\vert \sigma _{n}^{\alpha }a\left( x\right) \right\vert \leq \frac{%
c_{\alpha ,p}2^{M\left( 1/p-1\right) }2^{\alpha l+k}}{n^{\alpha }}.
\label{12}
\end{equation}

Let $x\in I_{M}^{k,M},\,0\leq k<M.$ Then from Lemma \ref{bt1} we have 
\begin{equation}
\left\vert \sigma _{n}^{\alpha }a\left( x\right) \right\vert \leq c_{\alpha
,p}2^{M\left( 1/p-1\right) +k}.  \label{12a}
\end{equation}

\bigskip By combining (\ref{2}), (\ref{12}) and (\ref{12a}) we obtain

\begin{equation*}
\int_{\overline{I_{M}}}\underset{n\in \mathbb{N}}{\sup }\left\vert \frac{%
\sigma _{n}^{\alpha }a\left( x\right) }{n^{1/p-1-\alpha }}\right\vert
^{p}d\mu \left( x\right) 
\end{equation*}%
\begin{equation*}
=\overset{M-2}{\underset{k=0}{\sum }}\overset{M-1}{\underset{l=k+1}{\sum }}%
\sum\limits_{x_{j}=0,j\in \{l+1,\dots ,M-1\}}^{1}\int_{I_{M}^{k,l}}\underset{%
n>2^{M}}{\sup }\left\vert \frac{\sigma _{n}^{\alpha }a\left( x\right) }{%
n^{1/p-1-\alpha }}\right\vert ^{p}d\mu \left( x\right) 
\end{equation*}%
\begin{equation*}
+\overset{M-1}{\underset{k=0}{\sum }}\int_{I_{M}^{k,M}}\underset{n>2^{M}}{%
\sup }\left\vert \frac{\sigma _{n}^{\alpha }a\left( x\right) }{%
n^{1/p-1-\alpha }}\right\vert ^{p}d\mu \left( x\right) 
\end{equation*}%
\begin{equation*}
\leq \frac{1}{2^{M\left( 1-\left( 1+\alpha \right) p\right) }}\overset{M-2}{%
\underset{k=0}{\sum }}\overset{M-1}{\underset{l=k+1}{\sum }}%
\sum\limits_{x_{j}=0,j\in \{l+1,\dots ,N-1\}}^{1}\int_{I_{M}^{k,l}}\underset{%
n>2^{M}}{\sup }\left\vert \sigma _{n}^{\alpha }a\left( x\right) \right\vert
^{p}d\mu \left( x\right) 
\end{equation*}%
\begin{equation*}
+\frac{1}{2^{M\left( 1-\left( 1+\alpha \right) p\right) }}\overset{M-1}{%
\underset{k=0}{\sum }}\int_{I_{M}^{k,M}}\underset{n>2^{M}}{\sup }\left\vert
\sigma _{n}^{\alpha }a\left( x\right) \right\vert ^{p}d\mu \left( x\right) 
\end{equation*}%
\begin{equation*}
\leq \frac{c_{\alpha ,p}}{2^{M\left( 1-\left( 1+\alpha \right) p\right) }}%
\overset{M-2}{\underset{k=0}{\sum }}\overset{M-1}{\underset{l=k+1}{\sum }}%
\frac{1}{2^{l}}\frac{2^{M\left( 1-p\right) }2^{\left( \alpha l+k\right) p}}{%
2^{\alpha pM}}
\end{equation*}%
\begin{equation*}
+\frac{c_{\alpha ,p}}{2^{M\left( 1-\left( 1+\alpha \right) p\right) }}\frac{1%
}{2^{M}}\overset{M-1}{\underset{k=0}{\sum }}2^{M\left( 1-p\right) +pk}
\end{equation*}%
\begin{equation*}
\leq \frac{c_{\alpha ,p}2^{\alpha pM}}{2^{\alpha pM}}\overset{M-2}{\underset{%
k=0}{\sum }}2^{kp}\overset{M-1}{\underset{l=k+1}{\sum }}\frac{1}{2^{l\left(
1-\alpha p\right) }}
\end{equation*}%
\begin{equation*}
+\frac{c_{\alpha ,p}}{2^{M\left( 1-\left( 1+\alpha \right) p\right) }}%
\overset{M-1}{\underset{k=0}{\sum }}\frac{2^{pk}}{2^{pM}}\leq c_{\alpha
,p}<\infty .
\end{equation*}

\textbf{b) }It is easy to show that under condition (\ref{6}), there exists
a sequence of positive integers$\ \left\{ n_{k},\ k\in \mathbb{N}%
_{+}\right\} ,$ such that 
\begin{equation*}
\lim_{k\rightarrow \infty }\frac{\left( 2^{2n_{k}}+2\right) ^{1/p-1-\alpha }%
}{\varphi \left( 2^{2n_{k}}+2\right) }=\infty .
\end{equation*}

Let 
\begin{equation*}
f_{n_{k}}=D_{2^{2n_{k}+1}}-D_{2^{_{2n_{k}}}}.
\end{equation*}

It is evident 
\begin{equation*}
\widehat{f}_{n_{k}}\left( i\right) =\left\{ 
\begin{array}{l}
1,\text{ if }i=2^{2n_{k}},\ldots ,2^{2n_{k}+1}-1, \\ 
0,\text{otherwise}.%
\end{array}%
\right.
\end{equation*}
Then we can write 
\begin{equation}
S_{i}f_{n_{k}}=\left\{ 
\begin{array}{ll}
D_{i}-D_{2^{2n_{k}}}, & \text{ if }i=2^{2n_{k}}+1,\ldots ,2^{2n_{k}+1}-1, \\ 
f_{n_{k}}, & \text{ if }i\geq 2^{2n_{k}+1}, \\ 
0, & \text{ otherwise.}%
\end{array}%
\right.  \label{33}
\end{equation}

From (\ref{1dn}) we get 
\begin{equation}
\left\Vert f_{n_{k}}\right\Vert _{H_{p}}=\left\Vert f_{n_{k}}^{\ast
}\right\Vert _{p}=\left\Vert D_{2^{2n_{k}+1}}-D_{2^{2n_{k}}}\right\Vert
_{p}\leq c2^{2n_{k}\left( 1-1/p\right) }.  \label{34}
\end{equation}

Since $A_{0}^{\alpha -1}=1,$ by (\ref{33}) we can write 
\begin{equation*}
\frac{\left\vert \sigma _{2^{2n_{k}}+1}^{\alpha }f_{n_{k}}\right\vert }{%
\varphi \left( 2^{2n_{k}}+1\right) }
\end{equation*}%
\begin{equation*}
=\frac{1}{\varphi \left( 2^{2n_{k}}+1\right) A_{2^{2n_{k}}+1}^{\alpha }}%
\left\vert \overset{2^{2n_{k}}+1}{\underset{j=1}{\sum }}A_{2^{2n_{k}}+1-j}^{%
\alpha -1}S_{j}f_{n_{k}}\right\vert
\end{equation*}%
\begin{equation*}
=\frac{1}{\varphi \left( 2^{2n_{k}}+1\right) A_{2^{2n_{k}}+1}^{\alpha }}%
\left\vert \overset{2^{2n_{k}}+1}{\underset{j=2^{2n_{k}}+1}{\sum }}%
A_{2^{2n_{k}}+1-j}^{\alpha -1}S_{j}f_{n_{k}}\right\vert
\end{equation*}%
\begin{equation*}
=\frac{1}{\varphi \left( 2^{2n_{k}}+1\right) A_{2^{2n_{k}}+1}^{\alpha }}%
\left\vert A_{0}^{\alpha -1}\left( D_{2^{2n_{k}}+1}-D_{2^{2n_{k}}}\right)
\right\vert
\end{equation*}%
\begin{equation*}
=\frac{1}{\varphi \left( 2^{2n_{k}}+1\right) A_{2^{2n_{k}}+1}^{\alpha }}%
\left\vert A_{0}^{\alpha -1}w_{2^{2n_{k}}}\right\vert
\end{equation*}%
\begin{equation*}
\geq \frac{c}{\varphi \left( 2^{2n_{k}}+1\right) \left( 2^{2n_{k}}+1\right)
^{\alpha }}.
\end{equation*}

From (\ref{34}) we have 
\begin{equation*}
\frac{c/\left( \varphi \left( 2^{2n_{k}}+1\right) \left( 2^{2n_{k}}+1\right)
^{\alpha }\right) \mu \left\{ x:\left\vert \overset{\sim }{\sigma }^{\alpha
,\ast }f\right\vert \geq c/\left( \varphi \left( 2^{2n_{k}}+1\right) \left(
2^{2n_{k}}+1\right) ^{\alpha }\right) \right\} ^{1/p}}{\left\Vert
f_{n_{k}}\right\Vert _{H_{p}}}
\end{equation*}%
\begin{equation*}
\geq \frac{c}{\varphi \left( 2^{2n_{k}}+1\right) \left( 2^{2n_{k}}+1\right)
^{\alpha }}\frac{1}{2^{2n_{k}\left( 1-1/p\right) }}\geq \frac{c\left(
2^{2n_{k}}+1\right) ^{1/p-1-\alpha }}{\varphi \left( 2^{2n_{k}}+1\right) }%
\rightarrow \infty ,\text{ as }k\rightarrow \infty .
\end{equation*}%
Theorem \ref{Theorem1} is proven.

\ 

\textbf{Proof of Theorem 2. }Suppose that 
\begin{equation*}
\overset{\infty }{\underset{m=1}{\sum }}\frac{\left\Vert \sigma _{m}^{\alpha
}F\right\Vert _{p}^{p}}{m^{2-\left( 1+\alpha \right) p}}\leq \left\Vert
F\right\Vert _{H_{p}}^{p}.
\end{equation*}
Then by using (\ref{5.1}) we have 
\begin{equation}
\overset{\infty }{\underset{m=1}{\sum }}\frac{\left\Vert \sigma _{m}^{\alpha
}F\right\Vert _{H_{p}}^{p}}{m^{2-\left( 1+\alpha \right) p}}  \label{5.3}
\end{equation}%
\begin{equation*}
=\overset{\infty }{\underset{m=1}{\sum }}\frac{\int_{G}\left\Vert \widetilde{%
\sigma _{m}^{\alpha }F^{\left( t\right) }}\right\Vert _{p}^{p}dt}{%
m^{2-\left( 1+\alpha \right) p}}\leq \int_{G}\overset{n}{\underset{m=1}{\sum 
}}\frac{\left\Vert \sigma _{m}^{\alpha }\widetilde{F^{\left( t\right) }}%
\right\Vert _{p}^{p}}{m^{2-\left( 1+\alpha \right) p}}dt
\end{equation*}%
\begin{equation*}
\leq \int_{G}\left\Vert \widetilde{F^{\left( t\right) }}\right\Vert
_{H_{p}}^{p}dt\sim \int_{G}\left\Vert F\right\Vert _{H_{p}}^{p}dt=\left\Vert
F\right\Vert _{H_{p}}^{p}.
\end{equation*}

According to Theorem \ref{W} and (\ref{5.3}) the proof of Theorem \ref%
{Theorem2} will be complete, if we show 
\begin{equation*}
\overset{\infty }{\underset{m=1}{\sum }}\frac{\left\Vert \sigma _{m}^{\alpha
}a\right\Vert _{p}^{p}}{m^{2-\left( 1+\alpha \right) p}}\leq c_{\alpha
}<\infty ,
\end{equation*}%
for every $p$-atom $a.$ Analogously to first part of Theorem 1 we can assume
that $n>2^{M}$ and $a$ be an arbitrary $p$-atom, with support$\ I,$ $\mu
\left( I\right) =2^{-M}$ \ and $I=I_{M}.$

Let $x\in I_{M}.$ Since $\sigma _{n}$ is bounded from $L_{\infty }$ to $%
L_{\infty }$ (the boundedness follows from (\ref{4})) and $\left\Vert
a\right\Vert _{\infty }\leq c2^{M/p}$ we obtain 
\begin{equation*}
\int_{I_{M}}\left\vert \sigma _{m}^{\alpha }a\right\vert ^{p}d\mu \leq
\int_{I_{M}}\left\Vert K_{m}^{\alpha }\right\Vert _{1}^{p}\left\Vert
a\right\Vert _{\infty }^{p}d\mu
\end{equation*}%
\begin{equation*}
\leq c_{\alpha ,p}\int_{I_{M}}\left\Vert a\right\Vert _{\infty }^{p}d\mu
\leq c_{\alpha ,p}<\infty .
\end{equation*}%
Hence%
\begin{equation*}
\overset{\infty }{\underset{m=2^{M}+1}{\sum }}\frac{\int_{I_{M}}\left\vert
\sigma _{m}^{\alpha }a\right\vert ^{p}d\mu }{m^{2-\left( 1+\alpha \right) p}}
\end{equation*}%
\begin{equation*}
\leq c_{\alpha ,p}\overset{\infty }{\underset{m=2^{M}+1}{\sum }}\frac{1}{%
m^{2-\left( 1+\alpha \right) p}}
\end{equation*}%
\begin{equation*}
\leq \frac{c_{\alpha ,p}}{2^{M\left( 1-\left( 1+\alpha \right) p\right) }}%
\leq c_{\alpha ,p}<\infty .
\end{equation*}

By combining (\ref{2}), (\ref{12}) and (\ref{12a}) analogously to first part
of Theorem 1 we can write%
\begin{equation*}
\overset{\infty }{\underset{m=2^{M}+1}{\sum }}\frac{\int_{\overline{I_{M}}%
}\left\vert \sigma _{m}^{\alpha }a\right\vert ^{p}d\mu }{m^{2-\left(
1+\alpha \right) p}}
\end{equation*}%
\begin{equation*}
=\overset{\infty }{\underset{m=2^{M}+1}{\sum }}\left( \overset{M-2}{\underset%
{k=0}{\sum }}\overset{M-1}{\underset{l=k+1}{\sum }}\sum\limits_{x_{j}=0,j\in
\{l+1,\dots ,M-1\}}^{1}\frac{\int_{I_{M}^{k,l}}\left\vert \sigma
_{m}^{\alpha }a\right\vert ^{p}d\mu }{m^{2-\left( 1+\alpha \right) p}}+%
\overset{M-1}{\underset{k=0}{\sum }}\frac{\int_{I_{M}^{k,M}}\left\vert
\sigma _{m}^{\alpha }a\right\vert ^{p}d\mu }{m^{2-\left( 1+\alpha \right) p}}%
\right)
\end{equation*}%
\begin{equation*}
\leq \overset{\infty }{\underset{m=2^{M}+1}{\sum }}\left( \frac{c_{\alpha
,p}2^{M\left( 1-p\right) }}{m^{2-p}}\overset{M-2}{\underset{k=0}{\sum }}%
\overset{M-1}{\underset{l=k+1}{\sum }}\frac{2^{p\left( \alpha l+k\right) }}{%
2^{l}}+\frac{c_{\alpha ,p}2^{M\left( 1-p\right) }}{m^{2-\left( 1+\alpha
\right) p}}\overset{M-1}{\underset{k=0}{\sum }}\frac{2^{pk}}{2^{M}}\right)
\end{equation*}%
\begin{equation*}
<c_{\alpha ,p}2^{M\left( 1-p\right) }\overset{\infty }{\underset{m=2^{M}+1}{%
\sum }}\frac{1}{m^{2-p}}+c_{\alpha ,p}\overset{\infty }{\underset{m=2^{M}+1}{%
\sum }}\frac{1}{m^{2-\left( 1+\alpha \right) p}}\leq c_{\alpha ,p}<\infty,
\end{equation*}
which completes the proof of Theorem \ref{Theorem2}.

\end{document}